\documentclass{article}

\usepackage{amsmath}
\usepackage{amssymb}
\usepackage{graphicx}
\usepackage{wasysym}
\usepackage{url}
\usepackage{pdflscape}
\usepackage{amssymb}
\usepackage{authblk}

\begin{document}

\date{}

\title{Method of precision increase by averaging with application to numerical
differentiation}

\author[]{Andrej Liptaj\footnote{andrej.liptaj@savba.sk}}
\affil[]{Institute of Physics, Slovak academy of Sciences \authorcr D\'{u}bravsk\'{a} cesta 9, 845 11 Bratislava, Slovakia}

\maketitle

\begin{abstract}
If several independent algorithms for a computer-calculated quantity
exist, then one can expect their results (which differ because of
numerical errors) to follow approximately Gaussian distribution. The
mean of this distribution, interpreted as the value of the quantity
of interest, can be determined with better precision than what is
the precision provided by a single algorithm. Often, with lack of
enough independent algorithms, one can proceed differently: many practical
algorithms introduce a bias using a parameter, e.g. a small but finite
number to compute a limit or a large but finite number (cutoff) to
approximate infinity. One may vary such parameter of a single algorithm
and interpret the resulting numbers as generated by several algorithms.
A numerical evidence for the validity of this approach is shown for
differentiation.
\end{abstract}

\section{Introduction\label{sec:Introduction}}

Average understood as summation (divided by a constant) is under very
general assumptions subject to the central limit theorem. This can
be used in numerical computations for precision increase. Indeed,
if several independent algorithms for computing a quantity of interest
exist, each of them having certain numerical imprecision, one may
average the results and get a smaller error. Depending on circumstances, this procedure may be regarded as repeated unbiased independent measurements with
random errors and, for this scenario, the expected shrinking of
the error is
\begin{equation}
\sigma=\frac{\sqrt{\sum_{i=1}^{N}\sigma_{i}^{2}}}{N}\approx\frac{\sigma^{typical}}{\sqrt{N}},
\label{Eq_ErrRed} 
\end{equation}
where one expects the numerical errors of the various methods not to
be very different (all close to a typical value $\sigma^{typical}$).

The question of algorithm independence arises. Clearly, if each algorithm
from a given set leads to the same result then the algorithms are,
in the mathematical sense, fully correlated. However, for what concerns numerical errors, a different optics can be adopted: otherwise result-equivalent algorithms may differ a lot in the functions they use (and the corresponding register operations) which de-correlates their numerical uncertainties. It is reasonable to assume that the numerical
errors arise from technical details of the computer processing
and do not actually depend a lot on the global ``idea'' of a given
algorithm. Therefore one can reasonably assume that algorithms which
differ in ``technical'' sense provide practically uncorrelated numerical
errors of their results.

Unfortunately, in practice, one usually does not have many independent
methods to compute a given quantity. However what is often the case
is a biased parameter-dependent algorithm. The parameters allow to
approximate an ideal situation which is inaccessible via computers:
a small (but finite) $h$ can be used as a step in numerical differentiation
or integration, a large (but finite) $\Lambda$ lambda can be used
as a cutoff (approximating infinity). A natural
idea arises: one may use different (but reasonable) values for these
parameters, get in each case a result and apply the previous ideas.
Two issues can be addressed here: bias and error correlations.

Obviously, any numerical differentiation (as an example) with nonzero
step $h$ is biased and even if infinite-precision computers were
available the result would not be fully correct. The overall ``wrongness''
thus has two components: numerical errors and bias. Here an expectation
can be made: if the averaging helps to shrink the numerical errors
then one should tend to use, in the averaging approach, an algorithm
(its parameter values) with smaller bias. In other words: the averaging
cannot remove the bias but does remove numerical effects and so one
expects to get the most precise results for less biased algorithm
compared to the bias leading to the most precise results for a single (i.e.
non-averaged) algorithm.

The correlation of uncertainties of results from the parameter-changing
approach is something one can examine empirically. The case of numerical
differentiation studied in this text shows that their mutual independence
is large enough to provide substantial error reduction. 

In what follows, this text fully focuses on the numerical differentiation.
To honestly study the subject I will use several methods of numerical
differentiation. The corresponding issues will be
reviewed in Sec. \ref{sec:Numerical-differentiation-method}. Next,
in Sec. \ref{sec:Testing-and-results}, I will explain the heuristic
testing method and present its results. In the last sections I will
discuss different result-related observations and make summary and
conclusions.

\section{Numerical differentiation methods\label{sec:Numerical-differentiation-method}}

In this text I focus on the (ill-conditioned) numerical differentiation
of a general (differentiable) function, I will therefore ignore 
special recipes suited for special situations\footnote{E.g. the set of analytic functions and well-conditioned differentiation
based on the Cauchy theorem.}. To make sure that the error shrinking by averaging is not limited
to some specific algorithm, I propose to test it on three different
differentiation methods (with an appropriately chosen $h$):
\begin{itemize}
\item Averaged finite difference (AFD)
\begin{align*}
f'_{AFD}(x,h) & =\frac{1}{2}\left[\frac{f\left(x+h\right)-f\left(x\right)}{h}+\frac{f\left(x\right)-f\left(x-h\right)}{h}\right],\\
 & =\frac{f\left(x+h\right)-f\left(x-h\right)}{2h}.
\end{align*}
\item ``Five-point rule'' based on the Richardson extrapolation (RE, \cite{Richardson307,Richardson299})
\[
f'_{RE}(x,h)=\frac{f\left(x-2h\right)-8f\left(x-h\right)+8f\left(x+h\right)-f\left(x+2h\right)}{12h}.
\]
The implementation of the numerical differentiation is in many common
mathematical computer packages based on the Richardson extrapolation.
\item Lanczos differentiation by integration (LDI, \cite{Lanczos:1956:AA})
\[
f'_{LDI}(x,h)=\frac{3}{2h^{3}}\int_{x-h}^{x+h}\left(x-t\right)f\left(t\right)dt.
\]
To evaluate the integral I use, in my programs, the composite Boole's
rule \cite{Abramowitz:1974:HMF:1098650} with 16 equidistant points $\left\{ x_{i}\right\} _{i=1}^{i=16}$,
$x_{1}=x-h$, $x_{16}=x+h,$ $x_{i+1}-x_{i}=\triangle x$
\begin{align*}
\int_{x-h}^{x+h}f\left(t\right)dt & \approx2\:\triangle x\:\left(I_{1}+I_{2}+I_{3}+I_{4}\right)/45,\\
I_{1} & =7\left[f\left(x_{1}\right)+f\left(x_{16}\right)\right],\\
I_{2} & =32\sum_{i=1}^{i=7}f\left(x_{2i+1}\right),\\
I_{3} & =14\sum_{i=1}^{i=3}f\left(x_{4i}\right),\\
I_{4} & =12\sum_{i=0}^{i=3}f\left(x_{4i+2}\right).
\end{align*}
\end{itemize}
Let me index these methods by the letter $k$, $k\varepsilon\left\{ \text{AFD, RE, LDI}\right\} $.
I implement the averaging procedure in the straightforward way
\[
f_{k}^{\prime AV}(x)=\frac{1}{N}\sum_{i=1}^{N}f'_{k}(x,h_{i}),\:h_{i}\epsilon H,
\]
where the set $H$ is chosen in function of the $h$ used in the single
algorithm computation as follows:
\begin{itemize}
\item For AFD $H=\left[0.5h,\:1.5h\right]$ where two options are investigated
\begin{itemize}
\item $h_{i}$ is generated as a random number with uniform distribution
from the interval $H$ (noted $\text{AFD}_{MC}^{AV}$).
\item successive values of $h_{i}$ are generated such as to be equidistant
with $h_{1}=0.5h$ and $h_{N}=1.5h$ (noted $\text{AFD}_{ED}^{AV}$).
\end{itemize}
\item For RE and LDI $H=\left[0.5h,\:1.5h\right]$, where $h_{i}$ is generated
as a random number with uniform distribution from this interval (only
this option is investigated).
\end{itemize}
Use of random numbers seems to be a safer option if aiming uncorrelated
errors, yet regular division of the interval is tested also. For
testing purposes I use a program\footnote{The program can be, at least temporarily, downloaded from \url{http://www.dthph.sav.sk/fileadmin/user_upload/liptaj/differentiationAveraging.zip}  or requested from the author. I also greatly profited from the WxMaxima software.} written in the \emph{JAVA} programming language and double precision
variables.

\section{Testing and results\label{sec:Testing-and-results}}

\begin{figure}
\begin{centering}
\includegraphics[width=0.8\linewidth]{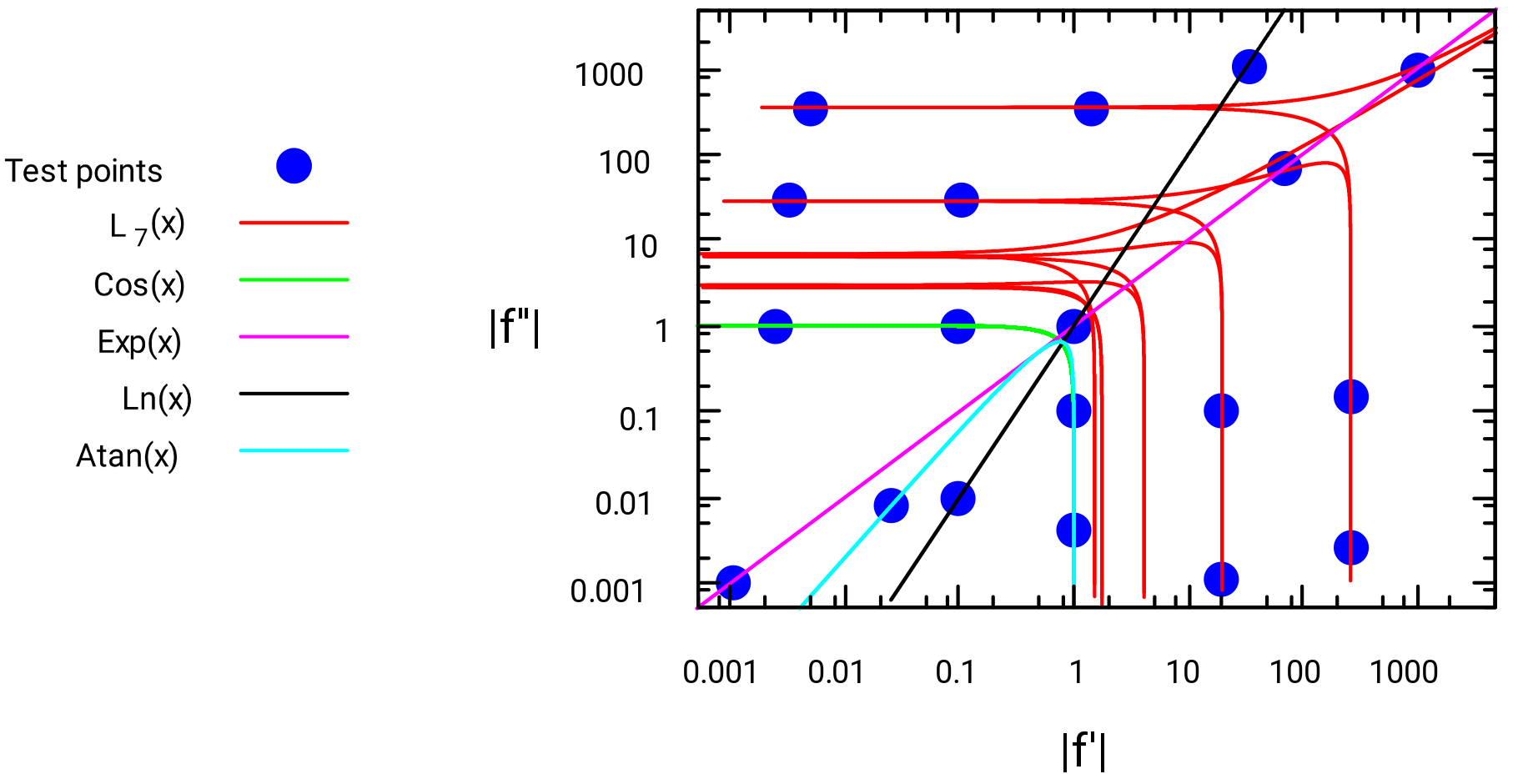}
\par\end{centering}
\caption{Studied cases depicted in $\left|f'\right|,\left|f''\right|$the plane.}
\label{Fig_derPlane}

\end{figure}

To study the behavior of the method in more details I make an effort to examine it depending on the first and second derivatives of the function and on the step size $h$. The first quantity directly correlates with what is being approximated ($f'$), the two others ($f'', h$) are often related to the expected precision of the approximation. I do the analysis by scanning 6 orders of magnitude for each ``dependence'' (its absolute value). For that purpose I choose 19 points in the $\left|f'\right|$, $\left|f''\right|$
plane, trying, in the logarithmic scale, to map it more or less uniformly.
To avoid any fine-tuning suspicions I choose to use the basic
elementary functions: $\cos\left(x\right)$, $\exp\left(x\right)$,
$\ln\left(x\right)$ and $\arctan\left(x\right)$. However, with this
choice, it is impossible to ``uniformly'' cover the $10^{-3}\apprle\left|f'\right|,\left|f''\right|\lesssim10^{3}$
region. Aiming this purpose, I add a suitable polynomial: the Laguerre
polynomial $L_{7}\left(x\right)$. Situation is summarized in Tab
\ref{Tab_19Cases} and in Fig. \ref{Fig_derPlane}. From now on I
will use the word ``case'' to refer to any of the 19 settings, each
of them characterized by a function $f$, its argument $x$ and the
absolute value of its first and second derivatives at $x$. I will
stick to the numbering presented in Tab. \ref{Tab_19Cases}.
\begin{table}
\centering{}%
\begin{tabular}{|c|c|c|c|c|}
\hline 
Case number & Function & $x$= & $\left|f'\left(x\right)\right|\approx$ & $\left|f''\left(x\right)\right|\approx$\tabularnewline
\hline 
\hline 
1 & $L_{7}\left(x\right)$ & $9.683$ & $19.88$ & $0.0011$\tabularnewline
\hline 
2 & $L_{7}\left(x\right)$ & $11.2345$ & $0.0031$ & $28.57$\tabularnewline
\hline 
3 & $L_{7}\left(x\right)$ & $15.83$ & $265.1$ & $0.1534$\tabularnewline
\hline 
4 & $L_{7}\left(x\right)$ & $17.65$ & $1.443$ & $358.1$\tabularnewline
\hline 
5 & $L_{7}\left(x\right)$ & $15.8285$ & $265.1$ & $0.0026$\tabularnewline
\hline 
6 & $L_{7}\left(x\right)$ & $17.64595$ & $0.0048$ & $356.8$\tabularnewline
\hline 
7 & $\exp\left(x\right)$ & $-6.9$ & $0.0010$ & $0.0010$\tabularnewline
\hline 
8 & $\ln\left(x\right)$ & $10$ & $0.1$ & $0.01$\tabularnewline
\hline 
9 & $\arctan\left(x\right)$ & $6.245$ & $0.0249$ & $0.0078$\tabularnewline
\hline 
10 & $\cos\left(x\right)$ & $1.47$ & $0.9949$ & $0.1006$\tabularnewline
\hline 
11 & $\cos\left(x\right)$ & $0.1$ & $0.0998$ & $0.9950$\tabularnewline
\hline 
12 & $\cos\left(x\right)$ & $0.0025$ & $0.0024$ & $0.9999$\tabularnewline
\hline 
13 & $\arctan\left(x\right)$ & $0.002$ & $0.9999$ & $0.0039$\tabularnewline
\hline 
14 & $\ln\left(x\right)$ & $0.03$ & $33.33$ & $1111.1$\tabularnewline
\hline 
15 & $\exp\left(x\right)$ & $6.9$ & $992.2$ & $992.2$\tabularnewline
\hline 
16 & $\ln\left(x\right)$ & $1.0$ & $1.0$ & $1.0$\tabularnewline
\hline 
17 & $L_{7}\left(x\right)$ & $9.67477$ & $19.88$ & $0.1000$\tabularnewline
\hline 
18 & $L_{7}\left(x\right)$ & $11.2311$ & $0.1001$ & $28.49$\tabularnewline
\hline 
19 & $\exp\left(x\right)$ & $4.25$ & $70.10$ & $70.10$\tabularnewline
\hline 
\end{tabular}\caption{Cases (points and functions) for which the averaging procedure was
tested.}
\label{Tab_19Cases}
\end{table}

The step size $h$ is changed from $h=10^{-3}$ to $h=10^{-8}$ in
geometrical progression with factor 10. One needs also to define the
size of the statistical sample. To profit most from the averaging
method a big number is suitable; I fix it to $N=10^{6}$. This choice
is driven also by practical considerations, i.e. the wish to keep
the computer processing time in reasonable limits ($\sim$ minutes).
The error is shown as absolute error
\begin{equation}
\left|f'_{\text{approximated}}-f'_{\text{true}}\right|,\label{eq:absRelErr}
\end{equation}
where for $f'_{\text{true}}$ the numerical value of the corresponding
(known) derivative function is taken.

To prevent long listings within the main text, I put the tables with
detailed results in Attachment. Each table corresponds to a single
case and is differential in the step size and used method. Here, I
average these tables (i.e. I average each cell over 19 cases), which
might be somewhat artificial but has more message-conveying power. 

\vspace{1cm}
{
\centering
\begin{tabular}{ l |  c c c}
\multicolumn{4}{c}{Case-averaged results}\\
$h$ & $10^{-3}$ & $10^{-4}$ & $10^{-5}$ \\ \hline \hline
$AFD$ & $ 2.0\times 10^{-3} $ & $2.0\times 10^{-5} $ & $\mathbf{1.8\times 10^{-7}} $\\
$AFD^{AV}_{MC}$ & $ 2.2\times 10^{-3} $ & $2.2\times 10^{-5} $ & $2.2\times 10^{-7} $\\
$AFD^{AV}_{ED}$ & $ 1.2\times 10^{-3} $ & $1.2\times 10^{-5} $ & $1.2\times 10^{-7} $ \\ \hline
$RE$ & $ \mathbf{6.9\times 10^{-9}} $ & $7.6\times 10^{-9} $ & $7.1\times 10^{-8} $ \\
$RE^{AV}$ & $ 1.1\times 10^{-8} $ & $\mathbf{6.0\times 10^{-12}} $ & $4.2\times 10^{-11} $ \\
 \hline $LDI$ & $ 1.2\times 10^{-3} $ & $\mathbf{1.4\times 10^{-5}} $ & $7.2\times 10^{-4} $ \\
$LDI^{AV}$ & $ 1.3\times 10^{-3} $ & $1.3\times 10^{-5} $ & $\mathbf{4.4\times 10^{-7}} $ \\ \hline
\end{tabular}

\vspace{1cm}

\begin{tabular}{ l |  c c c }
\multicolumn{4}{c}{Case-averaged results}\\
$h$ & $10^{-6}$ & $10^{-7}$ & $10^{-8}$ \\ \hline \hline
$AFD$ & $3.0\times 10^{-7} $ & $1.8\times 10^{-6} $ & $6.3\times 10^{-2} $ \\
$AFD^{AV}_{MC}$ & $\mathbf{2.3\times 10^{-9}} $ & $3.1\times 10^{-9} $ & $5.1\times 10^{-5} $ \\
$AFD^{AV}_{ED}$ &  $1.7\times 10^{-9} $ & $\mathbf{1.6\times 10^{-9}} $ & $1.0\times 10^{-6} $ \\
 \hline $RE$ &  $4.3\times 10^{-7} $ & $2.4\times 10^{-6} $ & $9.4\times 10^{-2} $ \\
$RE^{AV}$ &  $7.1\times 10^{-10} $ & $1.0\times 10^{-9} $ & $4.0\times 10^{-5} $ \\
 \hline $LDI$ & $7.2\times 10^{-2} $ & $5.5\times 10^{0} $ & $1.2\times 10^{7} $ \\
$LDI^{AV}$ &  $9.3\times 10^{-5} $ & $5.1\times 10^{-3} $ & $8.0\times 10^{3} $ \\
 \hline \end{tabular}

\vspace{1cm}
}
\section{Discussion}

Results confirm that the averaging method is very efficient in providing
precise numerical derivative and reducing related errors. The overall
error reduction (in absolute error) typically corresponds to two or
three orders of magnitude\footnote{One may notice that this error reduction roughly corresponds to what is predicted by the formula (\ref{Eq_ErrRed}). It gives additional hint in favor of the expected "modus operandi" of the averaging method.} (when comparing the most precise results).
Besides the obvious fact of reducing the error by increasing statistics\footnote{Non-averaged results can be seen as averaged results with statistics
equal to one.}, the assumptions concerning method functioning are further confirmed
by the behavior with respect to $h$: as predicted earlier (Sec. \ref{sec:Introduction})
the most precise results of the averaging method typically happen
for smaller step $h$ than is $h$ which corresponds to the most precise
result of the same, but non-averaged method. Rather numerous are situations
where $h$ remains the same, rare are exceptions where the behavior
is opposite ($RE$ method in case 12 and $AFD_{MC}^{AV}$ method in
case 16). Rather amazing are results of the $RE^{AV}$ method in cases
13 and 14 where, within the computer precision, exact results are
reconstructed.

When comparing $AFD_{MC}^{AV}$ and $AFD_{ED}^{AV}$ approaches, one
observes that their performances are rather equivalent. Yet, the ``equidistant''
method performs somewhat better which is little bit surprising: one
can imagine that a regular division could introduce some correlation
into the numbers to be averaged and thus slightly spoil the results.
One can speculate that this behavior could be related to what is know
from quasi-Monte Carlo methods: random numbers are often distributed
quite unevenly, i.e. the ``low-discrepancy'' of the equidistant
method may be the reason for it to win. It might certainly be a good
idea for further studies to use, within the averaging method,  low-discrepancy
sequences.

The results also show that the averaging method can be combined with
any of the three proposed ``standard'' methods, which points once
more to the general statistical aspects of the method. One may notice
that the ``standard'' methods differ quite not only in the definition
but also in the optimal step size $h$ (step where the maximal precision
is reached). The most precise of them is clearly the one based on
the Richardson extrapolation.

Finally, one needs to remark that for the $LDI$ method in cases 14 and 16
the averaging method fails. I cannot think about a solid explanation, it might by a random accident or it might be somehow related to the fact of $LDI$ being the least precise of the studied methods (or any other feature of this method).  At least in the case 14 the non-averaged result it atypically precise for this method, which might be interpreted as a "luck". In case 16 the difference between results is small, making the averaging  failure not to be so "pronounced". In any case I need to stress that, despite these two observations, the averaging method works in general very well also for the $LDI$ algorithm.

\section{Summary and conclusion}

In this text I made a numerical study of the averaging method applied
to the numerical differentiation. A rigorous approach to the whole
idea would require a rigorous treatment of the floating-point arithmetic
in computer registers. If possible, such an approach would certainly
be very tedious with many assumptions and special cases. I believe
a numerical evidence is strong enough to make claims about the method
and its mechanism. The method is efficient and provides an important
precision increase. It is very general and robust because of its statistical
character. It should be used in situations where precision is the
priority, its main drawback, slowness, makes it not suitable for quick
computations. When combined with a high precision ``standard'' method,
the averaging method is, to my knowledge, the most precise numerical
differentiation method at the market today. 

\section*{Acknowledgment}

This work was partly supported by the Slovak Grant Agency for Sciences
VEGA, grant No. 2/0153/17 and by the Slovak Research and Development
Agency APVV, grant No. APVV-0463-12.

\bibliography{precisionNumDiff.bib}{}
\bibliographystyle{ieeetr}

\appendix

\section*{Appendix}

The following tables give detailed results for cases mentioned in
 Tab \ref{Tab_19Cases}. In each table the step $h$ is varied
from $10^{-3}$ to $10^{-8}$ in columns, in rows different methods
are presented (notation from Sec. \ref{sec:Numerical-differentiation-method}
is used). Individual cells contain absolute error (formula \ref{eq:absRelErr}), the most precise of them is,
for each method, shown in bold characters.

\begin{landscape}

\centering

\begin{tabular}{ l |  c c c c c c}
 \multicolumn{7}{c}{Case number: 1 }\\
$h$ & $10^{-3}$ & $10^{-4}$ & $10^{-5}$ & $10^{-6}$ & $10^{-7}$ & $10^{-8}$ \\ \hline \hline
$AFD$ & $ 2.1\times 10^{-4} $ & $ 2.1\times 10^{-6} $ & $ \mathbf{4.9\times 10^{-9}} $ & $ 2.1\times 10^{-7} $ & $ 1.2\times 10^{-6} $ & $ 9.5\times 10^{-5}$ \\
$AFD^{AV}_{MC}$ & $ 2.2\times 10^{-4} $ & $ 2.2\times 10^{-6} $ & $ 2.2\times 10^{-8} $ & $ \mathbf{2.3\times 10^{-10}} $ & $ 1.8\times 10^{-9} $ & $ 6.4\times 10^{-8}$ \\
$AFD^{AV}_{ED}$ & $ 1.2\times 10^{-4} $ & $ 1.2\times 10^{-6} $ & $ 1.2\times 10^{-8} $ & $ \mathbf{3.9\times 10^{-10}} $ & $ 5.2\times 10^{-10} $ & $ 2.1\times 10^{-9}$ \\
 \hline $RE$ & $ \mathbf{1.5\times 10^{-10}} $ & $ 4.4\times 10^{-10} $ & $ 2.3\times 10^{-8} $ & $ 2.9\times 10^{-7} $ & $ 1.9\times 10^{-6} $ & $ 1.5\times 10^{-4}$ \\
$RE^{AV}$ & $ 5.0\times 10^{-11} $ & $ \mathbf{1.7\times 10^{-12}} $ & $ 1.2\times 10^{-11} $ & $ 2.9\times 10^{-10} $ & $ 7.8\times 10^{-11} $ & $ 9.4\times 10^{-8}$ \\
 \hline $LDI$ & $ 1.2\times 10^{-4} $ & $ \mathbf{8.4\times 10^{-7}} $ & $ 2.6\times 10^{-5} $ & $ 2.6\times 10^{-3} $ & $ 1.3\times 10^{-1} $ & $ 2.6\times 10^{4}$ \\
$LDI^{AV}$ & $ 1.3\times 10^{-4} $ & $ 1.3\times 10^{-6} $ & $ \mathbf{2.8\times 10^{-8}} $ & $ 8.4\times 10^{-6} $ & $ 5.5\times 10^{-4} $ & $ 1.1\times 10^{1}$ \\
 \hline \end{tabular}
\vspace{1cm}

\begin{tabular}{ l |  c c c c c c}
 \multicolumn{7}{c}{Case number: 2 }\\
$h$ & $10^{-3}$ & $10^{-4}$ & $10^{-5}$ & $10^{-6}$ & $10^{-7}$ & $10^{-8}$ \\ \hline \hline
$AFD$ & $ 3.9\times 10^{-4} $ & $ 3.9\times 10^{-6} $ & $ \mathbf{1.7\times 10^{-8}} $ & $ 3.4\times 10^{-7} $ & $ 1.8\times 10^{-6} $ & $ 4.8\times 10^{-4}$ \\
$AFD^{AV}_{MC}$ & $ 4.2\times 10^{-4} $ & $ 4.2\times 10^{-6} $ & $ 4.2\times 10^{-8} $ & $ \mathbf{4.5\times 10^{-10}} $ & $ 2.1\times 10^{-9} $ & $ 4.7\times 10^{-7}$ \\
$AFD^{AV}_{ED}$ & $ 2.3\times 10^{-4} $ & $ 2.3\times 10^{-6} $ & $ 2.3\times 10^{-8} $ & $ \mathbf{4.0\times 10^{-10}} $ & $ 6.0\times 10^{-10} $ & $ 1.8\times 10^{-9}$ \\
 \hline $RE$ & $ 2.2\times 10^{-9} $ & $ \mathbf{4.4\times 10^{-10}} $ & $ 3.0\times 10^{-8} $ & $ 4.1\times 10^{-7} $ & $ 2.2\times 10^{-6} $ & $ 7.2\times 10^{-4}$ \\
$RE^{AV}$ & $ 2.8\times 10^{-9} $ & $ \mathbf{7.0\times 10^{-12}} $ & $ 9.5\times 10^{-12} $ & $ 4.8\times 10^{-10} $ & $ 2.8\times 10^{-9} $ & $ 2.7\times 10^{-7}$ \\
 \hline $LDI$ & $ 2.3\times 10^{-4} $ & $ \mathbf{1.6\times 10^{-6}} $ & $ 4.9\times 10^{-5} $ & $ 4.9\times 10^{-3} $ & $ 2.4\times 10^{-1} $ & $ 1.2\times 10^{5}$ \\
$LDI^{AV}$ & $ 2.5\times 10^{-4} $ & $ 2.5\times 10^{-6} $ & $ \mathbf{1.4\times 10^{-7}} $ & $ 5.5\times 10^{-6} $ & $ 4.9\times 10^{-4} $ & $ 2.6\times 10^{0}$ \\
 \hline \end{tabular}
\vspace{1cm}

\begin{tabular}{ l |  c c c c c c}
 \multicolumn{7}{c}{Case number: 3 }\\
$h$ & $10^{-3}$ & $10^{-4}$ & $10^{-5}$ & $10^{-6}$ & $10^{-7}$ & $10^{-8}$ \\ \hline \hline
$AFD$ & $ 1.7\times 10^{-3} $ & $ 1.7\times 10^{-5} $ & $ \mathbf{9.6\times 10^{-8}} $ & $ 2.5\times 10^{-7} $ & $ 1.1\times 10^{-6} $ & $ 4.6\times 10^{-2}$ \\
$AFD^{AV}_{MC}$ & $ 1.8\times 10^{-3} $ & $ 1.8\times 10^{-5} $ & $ 1.8\times 10^{-7} $ & $ \mathbf{1.3\times 10^{-9}} $ & $ 2.3\times 10^{-8} $ & $ 2.2\times 10^{-5}$ \\
$AFD^{AV}_{ED}$ & $ 9.8\times 10^{-4} $ & $ 9.8\times 10^{-6} $ & $ 9.8\times 10^{-8} $ & $ \mathbf{5.5\times 10^{-9}} $ & $ 1.4\times 10^{-8} $ & $ 7.4\times 10^{-7}$ \\
 \hline $RE$ & $ 1.1\times 10^{-8} $ & $ \mathbf{6.3\times 10^{-9}} $ & $ 1.0\times 10^{-7} $ & $ 4.8\times 10^{-7} $ & $ 5.1\times 10^{-7} $ & $ 6.9\times 10^{-2}$ \\
$RE^{AV}$ & $ 1.8\times 10^{-8} $ & $ \mathbf{2.0\times 10^{-11}} $ & $ 1.2\times 10^{-10} $ & $ 4.1\times 10^{-9} $ & $ 1.9\times 10^{-9} $ & $ 6.4\times 10^{-5}$ \\
 \hline $LDI$ & $ 1.0\times 10^{-3} $ & $ \mathbf{1.1\times 10^{-6}} $ & $ 6.0\times 10^{-4} $ & $ 6.0\times 10^{-2} $ & $ 3.0\times 10^{0} $ & $ 1.2\times 10^{7}$ \\
$LDI^{AV}$ & $ 1.1\times 10^{-3} $ & $ 1.1\times 10^{-5} $ & $ \mathbf{6.2\times 10^{-7}} $ & $ 4.6\times 10^{-5} $ & $ 6.5\times 10^{-3} $ & $ 1.0\times 10^{4}$ \\
 \hline \end{tabular}
\vspace{1cm}

\begin{tabular}{ l |  c c c c c c}
 \multicolumn{7}{c}{Case number: 4 }\\
$h$ & $10^{-3}$ & $10^{-4}$ & $10^{-5}$ & $10^{-6}$ & $10^{-7}$ & $10^{-8}$ \\ \hline \hline
$AFD$ & $ 5.3\times 10^{-3} $ & $ 5.3\times 10^{-5} $ & $ \mathbf{1.1\times 10^{-7}} $ & $ 2.6\times 10^{-6} $ & $ 1.9\times 10^{-5} $ & $ 5.4\times 10^{-1}$ \\
$AFD^{AV}_{MC}$ & $ 5.7\times 10^{-3} $ & $ 5.7\times 10^{-5} $ & $ 5.7\times 10^{-7} $ & $ \mathbf{4.9\times 10^{-9}} $ & $ 1.7\times 10^{-8} $ & $ 6.3\times 10^{-4}$ \\
$AFD^{AV}_{ED}$ & $ 3.1\times 10^{-3} $ & $ 3.1\times 10^{-5} $ & $ 3.1\times 10^{-7} $ & $ \mathbf{2.3\times 10^{-10}} $ & $ 9.7\times 10^{-10} $ & $ 9.0\times 10^{-6}$ \\
 \hline $RE$ & $ \mathbf{2.1\times 10^{-8}} $ & $ 8.7\times 10^{-8} $ & $ 6.2\times 10^{-7} $ & $ 3.8\times 10^{-6} $ & $ 2.7\times 10^{-5} $ & $ 8.0\times 10^{-1}$ \\
$RE^{AV}$ & $ 2.7\times 10^{-8} $ & $ \mathbf{1.8\times 10^{-11}} $ & $ 3.4\times 10^{-10} $ & $ 1.7\times 10^{-9} $ & $ 4.0\times 10^{-9} $ & $ 5.7\times 10^{-5}$ \\
 \hline $LDI$ & $ 3.2\times 10^{-3} $ & $ \mathbf{2.2\times 10^{-5}} $ & $ 9.6\times 10^{-4} $ & $ 9.6\times 10^{-2} $ & $ 9.6\times 10^{0} $ & $ 9.9\times 10^{7}$ \\
$LDI^{AV}$ & $ 3.4\times 10^{-3} $ & $ 3.4\times 10^{-5} $ & $ \mathbf{1.7\times 10^{-6}} $ & $ 7.1\times 10^{-4} $ & $ 1.9\times 10^{-2} $ & $ 6.5\times 10^{4}$ \\
 \hline \end{tabular}
\vspace{1cm}

\begin{tabular}{ l |  c c c c c c}
 \multicolumn{7}{c}{Case number: 5 }\\
$h$ & $10^{-3}$ & $10^{-4}$ & $10^{-5}$ & $10^{-6}$ & $10^{-7}$ & $10^{-8}$ \\ \hline \hline
$AFD$ & $ 1.7\times 10^{-3} $ & $ 1.7\times 10^{-5} $ & $ \mathbf{5.2\times 10^{-8}} $ & $ 2.1\times 10^{-7} $ & $ 5.8\times 10^{-6} $ & $ 4.7\times 10^{-2}$ \\
$AFD^{AV}_{MC}$ & $ 1.8\times 10^{-3} $ & $ 1.8\times 10^{-5} $ & $ 1.8\times 10^{-7} $ & $ \mathbf{2.2\times 10^{-9}} $ & $ 6.7\times 10^{-9} $ & $ 6.0\times 10^{-5}$ \\
$AFD^{AV}_{ED}$ & $ 9.8\times 10^{-4} $ & $ 9.8\times 10^{-6} $ & $ 9.8\times 10^{-8} $ & $ \mathbf{2.4\times 10^{-9}} $ & $ 5.7\times 10^{-9} $ & $ 1.4\times 10^{-7}$ \\
 \hline $RE$ & $ 9.2\times 10^{-9} $ & $ \mathbf{7.9\times 10^{-9}} $ & $ 1.6\times 10^{-7} $ & $ 4.5\times 10^{-7} $ & $ 5.9\times 10^{-6} $ & $ 6.9\times 10^{-2}$ \\
$RE^{AV}$ & $ 1.8\times 10^{-8} $ & $ \mathbf{1.9\times 10^{-11}} $ & $ 3.3\times 10^{-11} $ & $ 1.3\times 10^{-9} $ & $ 4.3\times 10^{-9} $ & $ 3.1\times 10^{-5}$ \\
 \hline $LDI$ & $ 1.0\times 10^{-3} $ & $ \mathbf{1.1\times 10^{-6}} $ & $ 6.0\times 10^{-4} $ & $ 6.0\times 10^{-2} $ & $ 3.0\times 10^{0} $ & $ 1.2\times 10^{7}$ \\
$LDI^{AV}$ & $ 1.1\times 10^{-3} $ & $ 1.1\times 10^{-5} $ & $ \mathbf{2.0\times 10^{-7}} $ & $ 3.5\times 10^{-5} $ & $ 2.2\times 10^{-3} $ & $ 1.0\times 10^{4}$ \\
 \hline \end{tabular}
\vspace{1cm}

\begin{tabular}{ l |  c c c c c c}
 \multicolumn{7}{c}{Case number: 6 }\\
$h$ & $10^{-3}$ & $10^{-4}$ & $10^{-5}$ & $10^{-6}$ & $10^{-7}$ & $10^{-8}$ \\ \hline \hline
$AFD$ & $ 5.3\times 10^{-3} $ & $ 5.3\times 10^{-5} $ & $ \mathbf{7.7\times 10^{-7}} $ & $ 1.5\times 10^{-6} $ & $ 1.6\times 10^{-6} $ & $ 5.6\times 10^{-1}$ \\
$AFD^{AV}_{MC}$ & $ 5.7\times 10^{-3} $ & $ 5.7\times 10^{-5} $ & $ 5.7\times 10^{-7} $ & $ 7.4\times 10^{-9} $ & $ \mathbf{2.2\times 10^{-9}} $ & $ 2.7\times 10^{-4}$ \\
$AFD^{AV}_{ED}$ & $ 3.1\times 10^{-3} $ & $ 3.1\times 10^{-5} $ & $ 3.1\times 10^{-7} $ & $ 8.9\times 10^{-9} $ & $ \mathbf{5.7\times 10^{-9}} $ & $ 8.9\times 10^{-6}$ \\
 \hline $RE$ & $ \mathbf{1.7\times 10^{-8}} $ & $ 3.3\times 10^{-8} $ & $ 3.3\times 10^{-7} $ & $ 1.9\times 10^{-6} $ & $ 2.0\times 10^{-6} $ & $ 8.4\times 10^{-1}$ \\
$RE^{AV}$ & $ 2.7\times 10^{-8} $ & $ \mathbf{3.0\times 10^{-11}} $ & $ 2.1\times 10^{-10} $ & $ 4.8\times 10^{-9} $ & $ 1.2\times 10^{-9} $ & $ 6.0\times 10^{-4}$ \\
 \hline $LDI$ & $ 3.2\times 10^{-3} $ & $ \mathbf{2.2\times 10^{-5}} $ & $ 9.6\times 10^{-4} $ & $ 9.6\times 10^{-2} $ & $ 9.6\times 10^{0} $ & $ 9.8\times 10^{7}$ \\
$LDI^{AV}$ & $ 3.4\times 10^{-3} $ & $ 3.4\times 10^{-5} $ & $ \mathbf{2.0\times 10^{-6}} $ & $ 2.8\times 10^{-4} $ & $ 6.2\times 10^{-2} $ & $ 6.7\times 10^{4}$ \\
 \hline \end{tabular}
\vspace{1cm}

\begin{tabular}{ l |  c c c c c c}
 \multicolumn{7}{c}{Case number: 7 }\\
$h$ & $10^{-3}$ & $10^{-4}$ & $10^{-5}$ & $10^{-6}$ & $10^{-7}$ & $10^{-8}$ \\ \hline \hline
$AFD$ & $ 2.1\times 10^{-2} $ & $ 2.1\times 10^{-4} $ & $ 2.1\times 10^{-6} $ & $ 2.5\times 10^{-7} $ & $ 1.1\times 10^{-7} $ & $ \mathbf{3.7\times 10^{-12}}$ \\
$AFD^{AV}_{MC}$ & $ 2.3\times 10^{-2} $ & $ 2.3\times 10^{-4} $ & $ 2.3\times 10^{-6} $ & $ 2.2\times 10^{-8} $ & $ 6.6\times 10^{-10} $ & $ \mathbf{5.0\times 10^{-15}}$ \\
$AFD^{AV}_{ED}$ & $ 1.2\times 10^{-2} $ & $ 1.2\times 10^{-4} $ & $ 1.2\times 10^{-6} $ & $ 1.2\times 10^{-8} $ & $ 9.6\times 10^{-10} $ & $ \mathbf{6.2\times 10^{-17}}$ \\
 \hline $RE$ & $ 3.1\times 10^{-8} $ & $ 4.0\times 10^{-9} $ & $ 2.4\times 10^{-8} $ & $ 4.0\times 10^{-7} $ & $ 8.0\times 10^{-7} $ & $ \mathbf{5.6\times 10^{-12}}$ \\
$RE^{AV}$ & $ 4.7\times 10^{-8} $ & $ 2.7\times 10^{-12} $ & $ 1.5\times 10^{-11} $ & $ 2.8\times 10^{-10} $ & $ 1.8\times 10^{-9} $ & $ \mathbf{4.2\times 10^{-15}}$ \\
 \hline $LDI$ & $ 1.3\times 10^{-2} $ & $ \mathbf{2.0\times 10^{-4}} $ & $ 1.0\times 10^{-2} $ & $ 1.0\times 10^{0} $ & $ 7.8\times 10^{1} $ & $ 2.7\times 10^{-4}$ \\
$LDI^{AV}$ & $ 1.4\times 10^{-2} $ & $ 1.4\times 10^{-4} $ & $ 3.5\times 10^{-6} $ & $ 6.5\times 10^{-4} $ & $ 4.0\times 10^{-3} $ & $ \mathbf{3.3\times 10^{-7}}$ \\
 \hline \end{tabular}
\vspace{1cm}

\begin{tabular}{ l |  c c c c c c}
 \multicolumn{7}{c}{Case number: 8 }\\
$h$ & $10^{-3}$ & $10^{-4}$ & $10^{-5}$ & $10^{-6}$ & $10^{-7}$ & $10^{-8}$ \\ \hline \hline
$AFD$ & $ 2.5\times 10^{-4} $ & $ 2.5\times 10^{-6} $ & $ \mathbf{4.8\times 10^{-8}} $ & $ 2.2\times 10^{-7} $ & $ 7.6\times 10^{-8} $ & $ 1.3\times 10^{-4}$ \\
$AFD^{AV}_{MC}$ & $ 2.7\times 10^{-4} $ & $ 2.7\times 10^{-6} $ & $ 2.7\times 10^{-8} $ & $ \mathbf{2.4\times 10^{-10}} $ & $ 1.3\times 10^{-9} $ & $ 7.7\times 10^{-8}$ \\
$AFD^{AV}_{ED}$ & $ 1.5\times 10^{-4} $ & $ 1.5\times 10^{-6} $ & $ 1.5\times 10^{-8} $ & $ \mathbf{7.2\times 10^{-11}} $ & $ 8.6\times 10^{-11} $ & $ 3.7\times 10^{-10}$ \\
 \hline $RE$ & $ \mathbf{2.7\times 10^{-10}} $ & $ 5.1\times 10^{-10} $ & $ 2.7\times 10^{-8} $ & $ 3.1\times 10^{-7} $ & $ 1.6\times 10^{-7} $ & $ 1.9\times 10^{-4}$ \\
$RE^{AV}$ & $ 5.0\times 10^{-10} $ & $ \mathbf{2.0\times 10^{-12}} $ & $ 1.2\times 10^{-11} $ & $ 2.6\times 10^{-11} $ & $ 4.2\times 10^{-10} $ & $ 1.4\times 10^{-7}$ \\
 \hline $LDI$ & $ 1.5\times 10^{-4} $ & $ \mathbf{1.0\times 10^{-6}} $ & $ 3.3\times 10^{-5} $ & $ 3.3\times 10^{-3} $ & $ 1.6\times 10^{-1} $ & $ 3.5\times 10^{4}$ \\
$LDI^{AV}$ & $ 1.6\times 10^{-4} $ & $ 1.6\times 10^{-6} $ & $ \mathbf{6.9\times 10^{-8}} $ & $ 8.7\times 10^{-6} $ & $ 7.1\times 10^{-4} $ & $ 1.7\times 10^{1}$ \\
 \hline \end{tabular}
\vspace{1cm}

\begin{tabular}{ l |  c c c c c c}
 \multicolumn{7}{c}{Case number: 9 }\\
$h$ & $10^{-3}$ & $10^{-4}$ & $10^{-5}$ & $10^{-6}$ & $10^{-7}$ & $10^{-8}$ \\ \hline \hline
$AFD$ & $ 8.6\times 10^{-5} $ & $ 8.6\times 10^{-7} $ & $ 8.5\times 10^{-9} $ & $ \mathbf{3.0\times 10^{-9}} $ & $ 2.7\times 10^{-8} $ & $ 1.3\times 10^{-6}$ \\
$AFD^{AV}_{MC}$ & $ 9.3\times 10^{-5} $ & $ 9.3\times 10^{-7} $ & $ 9.3\times 10^{-9} $ & $ \mathbf{7.0\times 10^{-11}} $ & $ 1.1\times 10^{-10} $ & $ 3.8\times 10^{-10}$ \\
$AFD^{AV}_{ED}$ & $ 5.0\times 10^{-5} $ & $ 5.0\times 10^{-7} $ & $ 5.0\times 10^{-9} $ & $ \mathbf{3.3\times 10^{-11}} $ & $ 4.3\times 10^{-11} $ & $ 1.3\times 10^{-10}$ \\
 \hline $RE$ & $ 1.1\times 10^{-9} $ & $ \mathbf{7.6\times 10^{-11}} $ & $ 2.0\times 10^{-10} $ & $ 6.2\times 10^{-9} $ & $ 4.2\times 10^{-8} $ & $ 1.7\times 10^{-6}$ \\
$RE^{AV}$ & $ 1.6\times 10^{-9} $ & $ 1.9\times 10^{-13} $ & $ \mathbf{1.9\times 10^{-13}} $ & $ 2.2\times 10^{-11} $ & $ 1.7\times 10^{-10} $ & $ 1.3\times 10^{-9}$ \\
 \hline $LDI$ & $ 5.2\times 10^{-5} $ & $ \mathbf{5.6\times 10^{-7}} $ & $ 5.5\times 10^{-6} $ & $ 5.5\times 10^{-4} $ & $ 4.1\times 10^{-2} $ & $ 1.4\times 10^{2}$ \\
$LDI^{AV}$ & $ 5.6\times 10^{-5} $ & $ 5.6\times 10^{-7} $ & $ \mathbf{7.1\times 10^{-9}} $ & $ 2.8\times 10^{-7} $ & $ 1.3\times 10^{-6} $ & $ 6.7\times 10^{-1}$ \\
 \hline \end{tabular}
\vspace{1cm}

\begin{tabular}{ l |  c c c c c c}
 \multicolumn{7}{c}{Case number: 10 }\\
$h$ & $10^{-3}$ & $10^{-4}$ & $10^{-5}$ & $10^{-6}$ & $10^{-7}$ & $10^{-8}$ \\ \hline \hline
$AFD$ & $ 4.5\times 10^{-5} $ & $ 4.5\times 10^{-7} $ & $ 4.5\times 10^{-9} $ & $ \mathbf{1.1\times 10^{-10}} $ & $ 6.7\times 10^{-10} $ & $ 1.4\times 10^{-10}$ \\
$AFD^{AV}_{MC}$ & $ 4.9\times 10^{-5} $ & $ 4.9\times 10^{-7} $ & $ 4.9\times 10^{-9} $ & $ 4.9\times 10^{-11} $ & $ \mathbf{2.6\times 10^{-12}} $ & $ 6.1\times 10^{-12}$ \\
$AFD^{AV}_{ED}$ & $ 2.7\times 10^{-5} $ & $ 2.7\times 10^{-7} $ & $ 2.7\times 10^{-9} $ & $ 2.7\times 10^{-11} $ & $ \mathbf{1.8\times 10^{-13}} $ & $ 2.9\times 10^{-13}$ \\
 \hline $RE$ & $ 3.9\times 10^{-9} $ & $ \mathbf{4.0\times 10^{-12}} $ & $ 4.1\times 10^{-12} $ & $ 2.1\times 10^{-10} $ & $ 1.4\times 10^{-9} $ & $ 6.0\times 10^{-10}$ \\
$RE^{AV}$ & $ 5.9\times 10^{-9} $ & $ 5.8\times 10^{-13} $ & $ \mathbf{1.7\times 10^{-14}} $ & $ 2.5\times 10^{-13} $ & $ 2.0\times 10^{-12} $ & $ 1.2\times 10^{-12}$ \\
 \hline $LDI$ & $ 2.7\times 10^{-5} $ & $ 2.7\times 10^{-7} $ & $ \mathbf{4.1\times 10^{-8}} $ & $ 7.6\times 10^{-6} $ & $ 1.1\times 10^{-3} $ & $ 8.7\times 10^{-1}$ \\
$LDI^{AV}$ & $ 3.0\times 10^{-5} $ & $ 3.0\times 10^{-7} $ & $ \mathbf{2.7\times 10^{-9}} $ & $ 2.6\times 10^{-8} $ & $ 1.7\times 10^{-6} $ & $ 1.0\times 10^{-3}$ \\
 \hline \end{tabular}
\vspace{1cm}

\begin{tabular}{ l |  c c c c c c}
 \multicolumn{7}{c}{Case number: 11 }\\
$h$ & $10^{-3}$ & $10^{-4}$ & $10^{-5}$ & $10^{-6}$ & $10^{-7}$ & $10^{-8}$ \\ \hline \hline
$AFD$ & $ 5.3\times 10^{-4} $ & $ 5.3\times 10^{-6} $ & $ 5.3\times 10^{-8} $ & $ 5.2\times 10^{-10} $ & $ 6.6\times 10^{-11} $ & $ \mathbf{3.9\times 10^{-11}}$ \\
$AFD^{AV}_{MC}$ & $ 5.7\times 10^{-4} $ & $ 5.7\times 10^{-6} $ & $ 5.7\times 10^{-8} $ & $ 5.7\times 10^{-10} $ & $ 5.7\times 10^{-12} $ & $ \mathbf{3.6\times 10^{-13}}$ \\
$AFD^{AV}_{ED}$ & $ 3.1\times 10^{-4} $ & $ 3.1\times 10^{-6} $ & $ 3.1\times 10^{-8} $ & $ 3.1\times 10^{-10} $ & $ 3.1\times 10^{-12} $ & $ \mathbf{5.1\times 10^{-15}}$ \\
 \hline $RE$ & $ 6.8\times 10^{-9} $ & $ 7.0\times 10^{-13} $ & $ \mathbf{5.8\times 10^{-13}} $ & $ 2.2\times 10^{-12} $ & $ 6.6\times 10^{-11} $ & $ 4.1\times 10^{-10}$ \\
$RE^{AV}$ & $ 1.0\times 10^{-8} $ & $ 1.0\times 10^{-12} $ & $ \mathbf{8.9\times 10^{-16}} $ & $ 3.6\times 10^{-15} $ & $ 4.5\times 10^{-14} $ & $ 1.4\times 10^{-13}$ \\
 \hline $LDI$ & $ 3.2\times 10^{-4} $ & $ 3.2\times 10^{-6} $ & $ \mathbf{3.3\times 10^{-8}} $ & $ 5.0\times 10^{-7} $ & $ 3.3\times 10^{-5} $ & $ 1.8\times 10^{-2}$ \\
$LDI^{AV}$ & $ 3.4\times 10^{-4} $ & $ 3.4\times 10^{-6} $ & $ 3.4\times 10^{-8} $ & $ \mathbf{2.2\times 10^{-10}} $ & $ 3.7\times 10^{-8} $ & $ 3.6\times 10^{-6}$ \\
 \hline \end{tabular}
\vspace{1cm}

\begin{tabular}{ l |  c c c c c c}
 \multicolumn{7}{c}{Case number: 12 }\\
$h$ & $10^{-3}$ & $10^{-4}$ & $10^{-5}$ & $10^{-6}$ & $10^{-7}$ & $10^{-8}$ \\ \hline \hline
$AFD$ & $ 5.8\times 10^{-4} $ & $ 5.8\times 10^{-6} $ & $ 5.8\times 10^{-8} $ & $ 5.8\times 10^{-10} $ & $ \mathbf{7.8\times 10^{-12}} $ & $ 4.3\times 10^{-10}$ \\
$AFD^{AV}_{MC}$ & $ 6.3\times 10^{-4} $ & $ 6.3\times 10^{-6} $ & $ 6.3\times 10^{-8} $ & $ 6.3\times 10^{-10} $ & $ 6.3\times 10^{-12} $ & $ \mathbf{2.8\times 10^{-13}}$ \\
$AFD^{AV}_{ED}$ & $ 3.4\times 10^{-4} $ & $ 3.4\times 10^{-6} $ & $ 3.4\times 10^{-8} $ & $ 3.4\times 10^{-10} $ & $ 3.4\times 10^{-12} $ & $ \mathbf{8.0\times 10^{-14}}$ \\
 \hline $RE$ & $ 7.0\times 10^{-9} $ & $ 7.6\times 10^{-13} $ & $ 5.4\times 10^{-13} $ & $ \mathbf{4.8\times 10^{-13}} $ & $ 4.5\times 10^{-11} $ & $ 6.1\times 10^{-10}$ \\
$RE^{AV}$ & $ 1.1\times 10^{-8} $ & $ 1.1\times 10^{-12} $ & $ \mathbf{1.8\times 10^{-15}} $ & $ 5.3\times 10^{-15} $ & $ 9.8\times 10^{-15} $ & $ 7.0\times 10^{-13}$ \\
 \hline $LDI$ & $ 3.5\times 10^{-4} $ & $ 3.5\times 10^{-6} $ & $ 3.5\times 10^{-8} $ & $ \mathbf{1.3\times 10^{-8}} $ & $ 2.6\times 10^{-6} $ & $ 1.3\times 10^{-4}$ \\
$LDI^{AV}$ & $ 3.8\times 10^{-4} $ & $ 3.8\times 10^{-6} $ & $ 3.8\times 10^{-8} $ & $ \mathbf{4.3\times 10^{-10}} $ & $ 4.4\times 10^{-9} $ & $ 3.8\times 10^{-7}$ \\
 \hline \end{tabular}
\vspace{1cm}

\begin{tabular}{ l |  c c c c c c}
 \multicolumn{7}{c}{Case number: 13 }\\
$h$ & $10^{-3}$ & $10^{-4}$ & $10^{-5}$ & $10^{-6}$ & $10^{-7}$ & $10^{-8}$ \\ \hline \hline
$AFD$ & $ 5.8\times 10^{-4} $ & $ 5.8\times 10^{-6} $ & $ 5.8\times 10^{-8} $ & $ 5.8\times 10^{-10} $ & $ 8.6\times 10^{-11} $ & $ \mathbf{3.7\times 10^{-11}}$ \\
$AFD^{AV}_{MC}$ & $ 6.3\times 10^{-4} $ & $ 6.3\times 10^{-6} $ & $ 6.3\times 10^{-8} $ & $ 6.3\times 10^{-10} $ & $ 6.4\times 10^{-12} $ & $ \mathbf{7.6\times 10^{-13}}$ \\
$AFD^{AV}_{ED}$ & $ 3.4\times 10^{-4} $ & $ 3.4\times 10^{-6} $ & $ 3.4\times 10^{-8} $ & $ 3.4\times 10^{-10} $ & $ 3.4\times 10^{-12} $ & $ \mathbf{4.3\times 10^{-14}}$ \\
 \hline $RE$ & $ 7.0\times 10^{-9} $ & $ 7.1\times 10^{-13} $ & $ \mathbf{2.7\times 10^{-13}} $ & $ 8.4\times 10^{-13} $ & $ 7.7\times 10^{-11} $ & $ 2.2\times 10^{-10}$ \\
$RE^{AV}$ & $ 1.1\times 10^{-8} $ & $ 1.1\times 10^{-12} $ & $ \mathbf{0.0\times 10^{0}} $ & $ 1.3\times 10^{-14} $ & $ 3.8\times 10^{-14} $ & $ 7.7\times 10^{-14}$ \\
 \hline $LDI$ & $ 3.5\times 10^{-4} $ & $ 3.5\times 10^{-6} $ & $ 3.5\times 10^{-8} $ & $ \mathbf{1.3\times 10^{-8}} $ & $ 2.6\times 10^{-6} $ & $ 1.3\times 10^{-4}$ \\
$LDI^{AV}$ & $ 3.8\times 10^{-4} $ & $ 3.8\times 10^{-6} $ & $ 3.8\times 10^{-8} $ & $ \mathbf{3.9\times 10^{-10}} $ & $ 2.1\times 10^{-9} $ & $ 7.2\times 10^{-7}$ \\
 \hline \end{tabular}
\vspace{1cm}

\begin{tabular}{ l |  c c c c c c}
 \multicolumn{7}{c}{Case number: 14 }\\
$h$ & $10^{-3}$ & $10^{-4}$ & $10^{-5}$ & $10^{-6}$ & $10^{-7}$ & $10^{-8}$ \\ \hline \hline
$AFD$ & $ 5.7\times 10^{-4} $ & $ 5.7\times 10^{-6} $ & $ 5.7\times 10^{-8} $ & $ 5.6\times 10^{-10} $ & $ \mathbf{1.9\times 10^{-11}} $ & $ 3.0\times 10^{-11}$ \\
$AFD^{AV}_{MC}$ & $ 6.1\times 10^{-4} $ & $ 6.1\times 10^{-6} $ & $ 6.1\times 10^{-8} $ & $ 6.1\times 10^{-10} $ & $ 6.2\times 10^{-12} $ & $ \mathbf{5.1\times 10^{-13}}$ \\
$AFD^{AV}_{ED}$ & $ 3.3\times 10^{-4} $ & $ 3.3\times 10^{-6} $ & $ 3.3\times 10^{-8} $ & $ 3.3\times 10^{-10} $ & $ 3.3\times 10^{-12} $ & $ \mathbf{4.9\times 10^{-14}}$ \\
 \hline $RE$ & $ 6.9\times 10^{-9} $ & $ 6.3\times 10^{-13} $ & $ \mathbf{2.5\times 10^{-14}} $ & $ 3.5\times 10^{-12} $ & $ 1.0\times 10^{-11} $ & $ 3.0\times 10^{-11}$ \\
$RE^{AV}$ & $ 1.0\times 10^{-8} $ & $ 1.0\times 10^{-12} $ & $ \mathbf{0.0\times 10^{0}} $ & $ 2.7\times 10^{-15} $ & $ 2.7\times 10^{-14} $ & $ 1.2\times 10^{-12}$ \\
 \hline $LDI$ & $ 3.4\times 10^{-4} $ & $ 3.4\times 10^{-6} $ & $ 3.6\times 10^{-8} $ & $ 2.5\times 10^{-7} $ & $ \mathbf{2.0\times 10^{-11}} $ & $ 1.1\times 10^{-3}$ \\
$LDI^{AV}$ & $ 3.7\times 10^{-4} $ & $ 3.7\times 10^{-6} $ & $ 3.7\times 10^{-8} $ & $ \mathbf{2.1\times 10^{-10}} $ & $ 6.2\times 10^{-9} $ & $ 4.8\times 10^{-6}$ \\
 \hline \end{tabular}
\vspace{1cm}

\begin{tabular}{ l |  c c c c c c}
 \multicolumn{7}{c}{Case number: 15 }\\
$h$ & $10^{-3}$ & $10^{-4}$ & $10^{-5}$ & $10^{-6}$ & $10^{-7}$ & $10^{-8}$ \\ \hline \hline
$AFD$ & $ 7.7\times 10^{-5} $ & $ 7.7\times 10^{-7} $ & $ 4.9\times 10^{-9} $ & $ \mathbf{3.1\times 10^{-9}} $ & $ 1.0\times 10^{-7} $ & $ 2.6\times 10^{-6}$ \\
$AFD^{AV}_{MC}$ & $ 8.3\times 10^{-5} $ & $ 8.3\times 10^{-7} $ & $ 8.3\times 10^{-9} $ & $ \mathbf{4.7\times 10^{-11}} $ & $ 7.6\times 10^{-11} $ & $ 3.6\times 10^{-9}$ \\
$AFD^{AV}_{ED}$ & $ 4.5\times 10^{-5} $ & $ 4.5\times 10^{-7} $ & $ 4.5\times 10^{-9} $ & $ 6.5\times 10^{-11} $ & $ 2.1\times 10^{-10} $ & $ \mathbf{6.5\times 10^{-11}}$ \\
 \hline $RE$ & $ 1.1\times 10^{-9} $ & $ \mathbf{2.2\times 10^{-10}} $ & $ 4.1\times 10^{-9} $ & $ 4.5\times 10^{-9} $ & $ 9.9\times 10^{-8} $ & $ 3.9\times 10^{-6}$ \\
$RE^{AV}$ & $ 1.8\times 10^{-9} $ & $ \mathbf{3.8\times 10^{-13}} $ & $ 1.6\times 10^{-12} $ & $ 9.6\times 10^{-12} $ & $ 4.8\times 10^{-10} $ & $ 1.4\times 10^{-9}$ \\
 \hline $LDI$ & $ 4.6\times 10^{-5} $ & $ \mathbf{5.1\times 10^{-7}} $ & $ 6.7\times 10^{-6} $ & $ 6.7\times 10^{-4} $ & $ 5.1\times 10^{-2} $ & $ 2.6\times 10^{2}$ \\
$LDI^{AV}$ & $ 5.0\times 10^{-5} $ & $ 5.0\times 10^{-7} $ & $ \mathbf{8.6\times 10^{-10}} $ & $ 3.7\times 10^{-7} $ & $ 2.8\times 10^{-5} $ & $ 1.1\times 10^{0}$ \\
 \hline \end{tabular}
\vspace{1cm}

\begin{tabular}{ l |  c c c c c c}
 \multicolumn{7}{c}{Case number: 16 }\\
$h$ & $10^{-3}$ & $10^{-4}$ & $10^{-5}$ & $10^{-6}$ & $10^{-7}$ & $10^{-8}$ \\ \hline \hline
$AFD$ & $ 1.6\times 10^{-4} $ & $ 1.6\times 10^{-6} $ & $ 1.6\times 10^{-8} $ & $ 1.3\times 10^{-10} $ & $ 1.2\times 10^{-10} $ & $ \mathbf{5.9\times 10^{-11}}$ \\
$AFD^{AV}_{MC}$ & $ 1.7\times 10^{-4} $ & $ 1.7\times 10^{-6} $ & $ 1.7\times 10^{-8} $ & $ 1.7\times 10^{-10} $ & $ \mathbf{6.7\times 10^{-13}} $ & $ 1.4\times 10^{-12}$ \\
$AFD^{AV}_{ED}$ & $ 9.1\times 10^{-5} $ & $ 9.1\times 10^{-7} $ & $ 9.1\times 10^{-9} $ & $ 9.1\times 10^{-11} $ & $ 1.1\times 10^{-12} $ & $ \mathbf{2.4\times 10^{-13}}$ \\
 \hline $RE$ & $ 4.8\times 10^{-9} $ & $ \mathbf{1.1\times 10^{-12}} $ & $ 5.8\times 10^{-12} $ & $ 3.8\times 10^{-11} $ & $ 3.0\times 10^{-10} $ & $ 1.2\times 10^{-9}$ \\
$RE^{AV}$ & $ 7.3\times 10^{-9} $ & $ 7.3\times 10^{-13} $ & $ \mathbf{1.8\times 10^{-14}} $ & $ 1.2\times 10^{-13} $ & $ 2.8\times 10^{-13} $ & $ 1.1\times 10^{-13}$ \\
 \hline $LDI$ & $ 9.4\times 10^{-5} $ & $ 9.4\times 10^{-7} $ & $ \mathbf{6.7\times 10^{-9}} $ & $ 3.4\times 10^{-7} $ & $ 2.0\times 10^{-5} $ & $ 4.4\times 10^{-2}$ \\
$LDI^{AV}$ & $ 1.0\times 10^{-4} $ & $ 1.0\times 10^{-6} $ & $ \mathbf{9.9\times 10^{-9}} $ & $ 2.3\times 10^{-8} $ & $ 2.4\times 10^{-6} $ & $ 1.6\times 10^{-2}$ \\
 \hline \end{tabular}
\vspace{1cm}

\begin{tabular}{ l |  c c c c c c}
 \multicolumn{7}{c}{Case number: 17 }\\
$h$ & $10^{-3}$ & $10^{-4}$ & $10^{-5}$ & $10^{-6}$ & $10^{-7}$ & $10^{-8}$ \\ \hline \hline
$AFD$ & $ 2.0\times 10^{-4} $ & $ 2.0\times 10^{-6} $ & $ 2.7\times 10^{-8} $ & $ \mathbf{1.8\times 10^{-8}} $ & $ 2.0\times 10^{-6} $ & $ 1.0\times 10^{-4}$ \\
$AFD^{AV}_{MC}$ & $ 2.2\times 10^{-4} $ & $ 2.2\times 10^{-6} $ & $ 2.2\times 10^{-8} $ & $ \mathbf{2.4\times 10^{-10}} $ & $ 1.2\times 10^{-9} $ & $ 6.9\times 10^{-10}$ \\
$AFD^{AV}_{ED}$ & $ 1.2\times 10^{-4} $ & $ 1.2\times 10^{-6} $ & $ 1.2\times 10^{-8} $ & $ \mathbf{1.3\times 10^{-10}} $ & $ 2.5\times 10^{-10} $ & $ 5.1\times 10^{-7}$ \\
 \hline $RE$ & $ \mathbf{1.8\times 10^{-10}} $ & $ 7.8\times 10^{-10} $ & $ 1.0\times 10^{-8} $ & $ 3.0\times 10^{-8} $ & $ 2.2\times 10^{-6} $ & $ 1.4\times 10^{-4}$ \\
$RE^{AV}$ & $ 3.9\times 10^{-11} $ & $ \mathbf{1.3\times 10^{-12}} $ & $ 2.2\times 10^{-11} $ & $ 2.8\times 10^{-10} $ & $ 2.2\times 10^{-9} $ & $ 6.7\times 10^{-8}$ \\
 \hline $LDI$ & $ 1.2\times 10^{-4} $ & $ \mathbf{8.3\times 10^{-7}} $ & $ 2.6\times 10^{-5} $ & $ 2.6\times 10^{-3} $ & $ 1.3\times 10^{-1} $ & $ 2.5\times 10^{4}$ \\
$LDI^{AV}$ & $ 1.3\times 10^{-4} $ & $ 1.3\times 10^{-6} $ & $ \mathbf{1.7\times 10^{-8}} $ & $ 2.7\times 10^{-7} $ & $ 1.2\times 10^{-4} $ & $ 3.6\times 10^{1}$ \\
 \hline \end{tabular}
\vspace{1cm}

\begin{tabular}{ l |  c c c c c c}
 \multicolumn{7}{c}{Case number: 18 }\\
$h$ & $10^{-3}$ & $10^{-4}$ & $10^{-5}$ & $10^{-6}$ & $10^{-7}$ & $10^{-8}$ \\ \hline \hline
$AFD$ & $ 3.9\times 10^{-4} $ & $ 3.9\times 10^{-6} $ & $ 2.9\times 10^{-8} $ & $ \mathbf{2.6\times 10^{-8}} $ & $ 1.4\times 10^{-6} $ & $ 4.8\times 10^{-4}$ \\
$AFD^{AV}_{MC}$ & $ 4.2\times 10^{-4} $ & $ 4.2\times 10^{-6} $ & $ 4.2\times 10^{-8} $ & $ \mathbf{8.5\times 10^{-10}} $ & $ 2.5\times 10^{-9} $ & $ 7.0\times 10^{-7}$ \\
$AFD^{AV}_{ED}$ & $ 2.3\times 10^{-4} $ & $ 2.3\times 10^{-6} $ & $ 2.3\times 10^{-8} $ & $ \mathbf{1.6\times 10^{-10}} $ & $ 2.4\times 10^{-9} $ & $ 9.8\times 10^{-9}$ \\
 \hline $RE$ & $ \mathbf{1.8\times 10^{-9}} $ & $ 3.2\times 10^{-9} $ & $ 1.7\times 10^{-8} $ & $ 6.3\times 10^{-8} $ & $ 1.9\times 10^{-6} $ & $ 6.8\times 10^{-4}$ \\
$RE^{AV}$ & $ 2.7\times 10^{-9} $ & $ \mathbf{7.1\times 10^{-12}} $ & $ 3.4\times 10^{-11} $ & $ 1.4\times 10^{-10} $ & $ 8.6\times 10^{-11} $ & $ 1.1\times 10^{-6}$ \\
 \hline $LDI$ & $ 2.3\times 10^{-4} $ & $ \mathbf{1.6\times 10^{-6}} $ & $ 4.9\times 10^{-5} $ & $ 4.9\times 10^{-3} $ & $ 2.4\times 10^{-1} $ & $ 1.2\times 10^{5}$ \\
$LDI^{AV}$ & $ 2.5\times 10^{-4} $ & $ 2.5\times 10^{-6} $ & $ \mathbf{3.7\times 10^{-8}} $ & $ 1.6\times 10^{-5} $ & $ 3.0\times 10^{-5} $ & $ 2.1\times 10^{2}$ \\
 \hline \end{tabular}
\vspace{1cm}

\begin{tabular}{ l |  c c c c c c}
 \multicolumn{7}{c}{Case number: 19 }\\
$h$ & $10^{-3}$ & $10^{-4}$ & $10^{-5}$ & $10^{-6}$ & $10^{-7}$ & $10^{-8}$ \\ \hline \hline
$AFD$ & $ 1.0\times 10^{-6} $ & $ 1.0\times 10^{-8} $ & $ \mathbf{7.6\times 10^{-10}} $ & $ 2.7\times 10^{-9} $ & $ 1.6\times 10^{-8} $ & $ 1.9\times 10^{-7}$ \\
$AFD^{AV}_{MC}$ & $ 1.1\times 10^{-6} $ & $ 1.1\times 10^{-8} $ & $ 1.1\times 10^{-10} $ & $ \mathbf{1.5\times 10^{-12}} $ & $ 5.7\times 10^{-11} $ & $ 2.2\times 10^{-10}$ \\
$AFD^{AV}_{ED}$ & $ 5.8\times 10^{-7} $ & $ 5.8\times 10^{-9} $ & $ 5.8\times 10^{-11} $ & $ 3.8\times 10^{-12} $ & $ 2.9\times 10^{-11} $ & $ \mathbf{3.0\times 10^{-12}}$ \\
 \hline $RE$ & $ 9.0\times 10^{-11} $ & $ \mathbf{6.3\times 10^{-11}} $ & $ 8.2\times 10^{-10} $ & $ 3.6\times 10^{-9} $ & $ 1.8\times 10^{-8} $ & $ 2.8\times 10^{-7}$ \\
$RE^{AV}$ & $ 1.4\times 10^{-10} $ & $ \mathbf{4.3\times 10^{-14}} $ & $ 8.0\times 10^{-13} $ & $ 1.3\times 10^{-12} $ & $ 4.2\times 10^{-11} $ & $ 2.2\times 10^{-10}$ \\
 \hline $LDI$ & $ 6.0\times 10^{-7} $ & $ \mathbf{1.7\times 10^{-8}} $ & $ 1.5\times 10^{-6} $ & $ 1.5\times 10^{-4} $ & $ 1.1\times 10^{-2} $ & $ 1.9\times 10^{1}$ \\
$LDI^{AV}$ & $ 6.5\times 10^{-7} $ & $ 6.5\times 10^{-9} $ & $ \mathbf{4.1\times 10^{-10}} $ & $ 4.1\times 10^{-8} $ & $ 1.3\times 10^{-5} $ & $ 8.0\times 10^{-2}$ \\
 \hline \end{tabular}

\end{landscape}

\end{document}